\documentclass{amsart}
\usepackage{amsmath,amssymb}
\usepackage{graphicx,subfigure,psfrag}

\newtheorem*{maintheorem}{Theorem}

\theoremstyle{definition}

\newtheorem*{ddefinition}{Definition}

\theoremstyle{remark}

\newtheorem*{rremark}{Remark}

\newcommand{\la}{\lambda}
\newcommand{\om}{\omega}

\def\RR{\mathbb{R}}

\newcommand{\cC}{{\mathcal C}}
\newcommand{\cD}{{\mathcal D}}

\newcommand{\cF}{{\mathcal F}}

\newcommand{\cH}{{\mathcal H}}
\newcommand{\cI}{{\mathcal I}}

\newcommand{\cR}{{\mathcal R}}

\newcommand{\pd}{\partial}
\newcommand\minus\backslash

\newcommand\lan\langle
\newcommand\ran\rangle

\DeclareMathOperator\Div{div}

\renewcommand\leq\leqslant
\renewcommand\geq\geqslant
\newlength{\intwidth}

 \DeclareMathOperator\curl{curl}

\newcommand\Vect{\mathfrak{X}^1_{\mathrm{ex}}}

\begin{document}

\title[Helicity is the only integral invariant]{Helicity is the only integral invariant\\ of volume-preserving transformations}

\author{Alberto Enciso}
\address{Instituto de Ciencias Matem\'aticas, Consejo Superior de
  Investigaciones Cient\'\i ficas, 28049 Madrid, Spain}
\email{aenciso@icmat.es, dperalta@icmat.es, fj.torres@icmat.es}

\author{Daniel Peralta-Salas}

\author{Francisco Torres de Lizaur}

%

\begin{abstract}
  We prove that any regular integral invariant of volume-preserving transformations is equivalent to the
  helicity. Specifically, given a functional~$\cI$ defined on
  exact divergence-free vector fields of class $C^1$ on a compact
  3-manifold that is associated with a well-behaved integral kernel,
  we prove that~$\cI$ is invariant under arbitrary volume-preserving
  diffeomorphisms if and only if it is a function of the helicity.
\end{abstract}
\maketitle

\noindent {\bf Significance statement:} Helicity is a remarkable conserved quantity that is fundamental to all the natural phenomena described by a vector field whose evolution is given by volume-preserving transformations. This is the case of the vorticity of an inviscid fluid flow or of the magnetic field of a conducting plasma. The topological nature of the helicity was unveiled by Moffatt, but its relevance goes well beyond that of being a new conservation law. Indeed, the helicity defines an integral invariant under any kind of volume-preserving diffeomorphisms. A well-known open problem is whether there exist any integral invariants other than the helicity. We answer this question by showing that, under some mild technical assumptions, the helicity is the only integral invariant.

\section{Introduction}

Incompressible inviscid fluids are modeled by the three-dimensional Euler equations,
which assert that the velocity field $u(x,t)$ of the fluid flow must
satisfy the system of differential equations
\[
\pd_t u+ (u\cdot \nabla ) u =-\nabla p\,,\qquad \Div u=0\,.
\]
Here the scalar function $p(x,t)$ is another unknown of the problem,
which physically corresponds to the pressure of the fluid. 

It is customary to introduce the vorticity $\om:=\curl u$ to simplify
the analysis of these equations, as it enables us to get rid of the
pressure function. In terms of the vorticity, the Euler equations read
as
\begin{equation}\label{pdtom}
\pd_t \om=[\om, u]\,,
\end{equation}
where $[\om,u]:=(\om\cdot \nabla ) u- (u\cdot \nabla) w$ is the
commutator of vector fields and $u$ can be written in terms of $\om$ using the
Biot--Savart law
\begin{equation}\label{BS}
u(x)=\curl^{-1}\om(x):=\frac1{4\pi}\int_{\RR^3}\frac{\om(y)\times(x-y)}{|x-y|^3}\, dy\,,
\end{equation}
at least when the space variable is assumed to take values in the
whole space~$\RR^3$.

The transport equation~\eqref{pdtom} was first derived by Helmholtz,
who showed that the meaning of this equation is that the vorticity at
time~$t$ is related to the vorticity at initial time $t_0$ via the
flow of the velocity field, provided
that the equation does not develop any singularities in the time interval
$[t_0,t]$. More precisely, if $\phi_{t,t_0}$ denotes the (time-dependent)
flow of the divergence-free field $u$, then the vorticity at time $t$
is given by the action of the push-forward of the volume-preserving diffeomorphism
$\phi_{t,t_0}$ on the initial vorticity:
\[
\om(\cdot, t)=(\phi_{t,t_0})_*\,\om(\cdot,t_0)\,.
\]

The phenomenon of the transport of vorticity gives rise to a new conservation law
of the three-dimensional Euler equations. Moffatt coined the term {\em helicity} for this conservation law in his influential paper~\cite{Mo69}, and exhibited its topological nature. Indeed, defining the
helicity of a divergence-free vector field~$w$ in $\RR^3$ as
\[
\cH(w):=\int_{\RR^3}w\cdot  \curl^{-1}w\, dx\,,
\]
it turns out that the helicity of the vorticity $\cH(\om(\cdot, t))$ is
a conserved quantity for the Euler equations. In fact, helicity is also conserved for the compressible Euler equations provided the fluid is barotropic (i.e. the pressure is a function of the density). 

It is well known that the relevance of the helicity goes well beyond
that of being a new (non-positive) conserved quantity for the Euler
equations. On the one hand, the helicity appears in other natural
phenomena that are also described by a divergence-free field whose
evolution is given by a time-dependent family of volume-preserving
diffeomorphisms~\cite{Mo14}. For instance, the case of magnetohydrodynamics (MHD), where
one is interested in the helicity of the magnetic field of a
conducting plasma, has attracted considerable attention. On the other
hand, it turns out that the helicity does not only correspond to a
conserved quantity for evolution equations such as Euler or MHD, but
in fact defines an integral invariant for vector fields under any kind
of volume-preserving diffeomorphisms~\cite{AK}. 

It is important to emphasize that conserved quantities of the Euler or MHD equations (e.g., the
kinetic energy and the momentum) are not, in
general, invariant under arbitrary volume-preserving diffeomorphisms, but they are invariant only under the very particular diffeomorphism defined by the flow of the velocity field of the fluid or conducting plasma. Perhaps the key feature of the helicity, which distinguishes it from other conserved quantities of Euler or MHD, is its invariance under any kind of volume-preserving transformations (in particular, it is invariant under the transport of the vorticity or the magnetic field by an arbitrary divergence-free vector field), so let us elaborate on this property.

Helicity is often analyzed in the context of a compact
3-dimensional manifold~$M$ without boundary, endowed with a Riemannian
metric. The simplest case would be that of the flat 3-torus, which
corresponds to fields on Euclidean space with periodic
boundary conditions. To define the helicity in a general compact
3-manifold, let us introduce some notation. We will denote by $\Vect$ the vector
space of exact divergence-free vector fields on $M$ of class $C^{1}$,
endowed with its natural $C^1$ norm. We recall that a divergence-free vector field~$w$
is {\em exact}\/ if its flux through any closed surface is zero (or,
equivalently, if there exists a vector field~$v$ such that $w=\curl
v$). This is a topological condition, and in particular when the first
homology group of the manifold is trivial (e.g., in the 3-sphere) every
divergence-free field is automatically exact.

As is well known, the reason to consider exact fields in this context is
that, on exact fields, the curl operator has a well defined inverse
$\curl^{-1}:\Vect\to\Vect$. The inverse of curl is a generalization to
compact 3-manifolds of the Biot--Savart operator~\eqref{BS}, and can
also be written in terms of a (matrix-valued) integral kernel $k(x,y)$ as
\begin{equation}\label{curl-1}
\curl^{-1} w (x)=\int_{M} k(x, y)\,w(y)\, dy\,,
\end{equation}
where $dy$ now stands for the Riemannian volume measure. 
Using this integral operator, one can define the helicity of a vector
field~$w$ on~$M$ as
$$
\cH(w):=\int_{M} w\cdot \curl^{-1} w \, dx\,.
$$
Here and in what follows the dot denotes the scalar product of
two vector fields defined by the Riemannian metric on $M$. The helicity is then invariant under volume-preserving
transformations, that is, $\cH(w)=\cH(\Phi_*w)$ for any
diffeomorphism~$\Phi$ of $M$ that preserves volume.

In view of the expression~\eqref{curl-1} for the inverse of the curl
operator, it is clear that the helicity is an {\em integral
  invariant}\/, meaning that it is given by the integral of a density
of the form
\[
\cH(w)=\int G(x,y,w(x),w(y))\, dx\, dy\,.
\]
Arnold and Khesin conjectured~\cite[Section I.9]{AK} that, in fact, the helicity is the only integral invariant, that is, there are no other invariants of the form
\begin{equation}\label{eq:intinv}
\cI(u):=\int G(x_{1},\dots,x_{n},
u(x_{1}),\dots,u(x_{n})) \ dx_{1}\cdots dx_{n}
\end{equation}
with $G$ a reasonably well-behaved function. Here all variables are assumed to be integrated over~$M$.

Our objective in this paper is to show, under some natural regularity
assumptions, that the helicity is indeed the only integral invariant
under volume-preserving diffeomorphisms. To this end, let us define a
{regular integral invariant} as follows:

\begin{ddefinition}\label{mainhypothesis}
Let $\cI:\Vect\to\RR$ be a $C^1$ functional. We say that $\cI$ is a
{\em regular integral invariant}\/ if:
\begin{enumerate}
\item It is invariant under volume-preserving
transformations, i.e., $\cI(w)=\cI(\Phi_*w)$ for any
diffeomorphism~$\Phi$ of $M$ that preserves volume.
\item At any point $w\in\Vect$, the (Fr\'echet) derivative of $\cI$ is an integral operator with continuous kernel, that is,
\[
(D\cI)_w(u)=\int_{M} K(w) \cdot u\,,
\]
for any $u\in\Vect$, where $K: \Vect\rightarrow\Vect $ is a continuous
map. 
\end{enumerate}
\end{ddefinition}

In the above definition and in what follows, we omit the Riemannian
volume measure under the integral sign when no
confusion can arise. Observe that any integral invariant of the form~\eqref{eq:intinv} is a
regular integral invariant provided that the function $G$ satisfies
some mild technical assumptions. In particular, the helicity is a regular integral invariant.

The following theorem, which is the main result of this paper, shows
 that the helicity is essentially the only regular integral invariant
 in the above sense. The proof of this result is presented in
 Section~\ref{S.main}, and is an extension to any closed
 3-manifold of a theorem of Kudryavtseva~\cite{Ku2}, who proved an
 analogous result for divergence-free vector fields on $3$-manifolds that are trivial bundles of a compact surface with boundary over the circle, which admit a cross section and are tangent to the boundary. Kudryavtseva's theorem is based on her work on the uniqueness of the Calabi invariant for area-preserving diffeomorphisms of the disk~\cite{Ku}. We observe that our main result does not imply the aforementioned theorem because we consider manifolds without boundary.

\begin{maintheorem}
Let $\cI$ be a regular integral invariant. Then $\cI$ is a function of the helicity, i.e., there exists a $C^{1}$ function $f:\RR\rightarrow\RR$ such that $\cI=f(\cH)$.
\end{maintheorem}

We would like to remark that this theorem does not exclude the
existence of other invariants of divergence-free vector fields under
volume-preserving diffeomorphisms that are not $C^1$ or whose
derivative is not an integral operator of the type described in
the definition above. For example, the KAM-type invariants
recently introduced in~\cite{KKP14} are in no way related to the
helicity, but they are not even continuous functionals on $\Vect$.

Other type of invariants that have attracted considerable attention
are the {asymptotic invariants} of divergence-free vector
fields~\cite{AR,Ghys,Kh03,B11,BM,Akhme,Komen}. These invariants are of
non-local nature because they are defined in terms of a knot invariant
(e.g., the linking number) and the flow of the vector field. In some
cases, it turns out that the asymptotic invariant can be expressed as
a regular integral invariant, as happens with the asymptotic linking
number for divergence-free vector fields~\cite{AR}, the asymptotic
signature~\cite{Ghys} and the asymptotic Vassiliev
invariants~\cite{BM,Komen} for ergodic divergence-free vector
fields. In these cases, the authors prove that the corresponding
asymptotic invariant is a function of the helicity, which is in
perfect agreement with our main theorem.  

The so-called {higher order helicities}~\cite{Be90,LS00,Komen09} are
also invariants under volume-preserving diffeomorphisms. However, they
are not defined for any divergence-free vector field, but just for
vector fields supported on a disjoint union of solid tori. This
property is, of course, not even continuous in $\Vect$, so these
functionals do not fall in the category of the regular integral invariants considered in this paper. 

Our main theorem is reminiscent of Serre's theorem~\cite{Se84} showing that any conserved quantity of the three-dimensional Euler equations that is the integral of a density depending on the velocity field and its first derivatives,
$$
\cI(u):=\int_{\RR^{3}}G(u(x,t), Du(x,t))\,dx\,,
$$
is a function of the energy, the momentum and the
helicity. From a technical point of view, 
the proof of our main theorem is totally different to the proof of
Serre's theorem, which is purely analytic, only holds in the Euclidean
space, and is based on integral identities that the density $G$ must
satisfy in order to define a conservation law of the Euler
equations. 

Even more importantly, from a conceptual standpoint it
should be emphasized that Serre's theorem applies to conserved
quantities of the Euler equations, while our theorem concerns the
existence of functionals that are invariant under any kind of volume-preserving diffeomorphisms, which is a much stronger requirement, as explained in a previous paragraph. In particular, the fact that the energy and the momentum are not functions of the helicity does not contradict our main theorem, because they are conserved by the evolution determined by the Euler
equations but they are not invariant under the flow of an
arbitrary divergence-free vector field. Accordingly, our theorem does not mean that there are no other integrals of motion of the Euler (or MHD) equations.  

It is worth noticing that one can construct well-behaved integral invariants of
Lagrangian type that are invariant under general volume-preserving
diffeomorphisms but which are not functions of the helicity. These
functional arise in a natural manner in the analysis of the Euler or
MHD equations especially when one considers integrable fields, that
is, fields whose integral curves are tangent to a family of invariant
surfaces. For example, one can define a partial helicity as the helicity integral taken over the region $\Omega$ bounded by an invariant surface of the field. In this context, if $f$ is any well-behaved function (e.g., a
smooth function supported on the region $\Omega$ covered by invariant surfaces)
which is assumed to be transported under the action of the
diffeomorphism group, the functional
\[
\cF(f,w):=\int_M f\, w\cdot \curl^{-1} w \, dx
\]
is invariant under volume-preserving diffeomorphisms (and it is not a
function of the helicity). The key point here is that the assumption
that $f$ is transformed in a Lagrangian way means that the action of
the volume-preserving diffeomorphism group is not the one considered
in this paper, which would be 
$$\Phi\cdot \cF(f,w):=\cF(f,\Phi_* w)\,,$$ 
but the one given by
\[
\Phi\cdot \cF(f,w):=\cF(f\circ \Phi^{-1},\Phi_* w)\,.
\]
In this sense, this new action is defined on functionals mapping a
function and a vector field (rather than just a vector field) to a number, so
it does not fall within the scope of our theorem. In the context of the partial helicity defined above, this action means that not only the vector field $w$, but also the function $f$ and the region $\Omega$ where it is supported, are transported by the fluid flow.

\section{Proof of the main theorem}
\label{S.main}

We divide the proof of the main theorem in five steps. The idea of the
proof, which is inspired by Kudryavtseva's work on the uniqueness of
the Calabi invariant~\cite{Ku}, is that the invariance of the functional $\cI$
under volume-preserving diffeomorphisms implies the existence of a
continuous first integral for each exact divergence-free vector
field. Since a generic vector field in~$\Vect$ is not integrable, we
conclude that the aforementioned first integral is a constant (that
depends on the field), which in turn implies that $\cI$ has the same
value for all vector fields in a connected component of the level sets
of the helicity. Since these level sets are path connected, the
theorem will follow.

\subsubsection*{Step 1: For each vector field $w\in\Vect$, either
  $\curl K(w)=fw$ on $M\backslash w^{-1}(0)$ for some 
  function~$f\in C^0(M\backslash w^{-1}(0))$ or the field~$w$ admits a nontrivial first
  integral (that is, $\nabla F\cdot w=0$ for some nonconstant function~$F\in C^1(M)$).}

We first notice that the flow $\phi_t$ of any divergence-free vector field $u$ is a $1$-parameter family of volume-preserving diffeomorphisms, so the functional $\cI$ must take the same values on $w$ and its push-forward $(\phi_t)_*w$, i.e.
$$
\cI((\phi_t)_*w)=\cI(w)
$$
for all $t\in\RR$. Taking derivatives with respect to $t$ in this equation and evaluating at $t=0$, we immediately get
\begin{equation}\label{eqpush}
0=\frac{d}{dt}\cI ((\phi_t)_*w)= (D\cI)_w([w, u])=\int_{M}K(w)\cdot[w, u]\,.
\end{equation}

The identity $[w, u]=\curl (u\times w)$ for divergence-free fields allows us to write the integral above as 
\begin{align*}
\int_{M} K(w)\cdot [w,u]&=\int_M K(w)\cdot \curl (u\times w)\\
&=\int_M \curl K(w)\cdot (u\times w)\\
&=\int_M u\cdot(w\times\curl K(w))
\end{align*}
where we have integrated by parts to obtain the second equality. Hence
Eq.~\eqref{eqpush} implies that for each pair of vector fields $u,w\in \Vect$ we have
$$
\int_M u\cdot(w\times\curl K(w))=0\,.
$$ 
It then follows that the vector field $w\times \curl K(w)$ is
$L^2$-orthogonal to all the divergence-free vector fields on $M$, and
hence the Hodge decomposition theorem implies that there exists a
$C^1$ function $F$ on $M$ such that $w\times \curl K(w)=\nabla
F$. Then $w\cdot \nabla F=0$, so $F$ is a first integral of $w$.  

In the case that $F$ is identically constant, we have that $w\times
\curl K(w)=0$, so $\curl K(w)$ is proportional to $w$ at any point of
$M$ where the latter does not vanish. Since $\curl K(w)$ is a
continuous vector field on $M$ because, by assumption, $K(w)\in \Vect$, it follows that there is a continuous function $f$ such that 
\begin{equation}\label{eqstep1}
\curl K(w)=fw
\end{equation}
in $M\backslash w^{-1}(0)$, as we wanted to prove.

\subsubsection*{Step 2: The function $f\in C^0(M\backslash w^{-1}(0))$ is a continuous first integral of $w$.} 

\begin{figure}[t]
  \centering
  \psfrag{D0}{$\mathcal D_0$}
\psfrag{D1}{$\mathcal D_1$}
  \psfrag{p0}{$p_0$}
  \psfrag{p1}{$p_1$}
  \psfrag{U}{$U$}
\includegraphics[scale=0.13,angle=0]{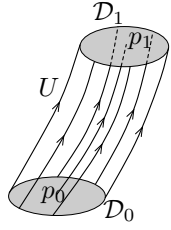}\hspace{0.1em}
\caption{A flow box for the vector field $w$.}\label{f1}
\end{figure}

The flow box theorem ensures that for any point in the complement of
the zero set $w^{-1}(0)$ there is a neighborhood $U$ and a
diffeomorphism $\Phi:U\to [0,1]\times D$ such that
$\Phi_*w=\partial_z$. Here $D:=\{x\in\RR^2:|x|\leq 1\}$ is the closed
unit 2-disk, and $[0,1]\times D$ is endowed with the natural Cartesian coordinates $x\in D$ and $z\in[0,1]$. Using the notation $\mathcal D_s:=\Phi^{-1}(\{s\}\times D)$ and $\mathcal S:=\Phi^{-1}([0,1]\times \partial D)$, it is obvious from the definition of the flow box that
$$
\partial U=\mathcal D_0\cup \mathcal D_1 \cup \mathcal S\,,
$$
and that the integral curves of $w$ are tangent to the cylinder
$\mathcal S$ and transverse to the disks $\mathcal D_0$ and $\mathcal
D_1$. This construction is depicted in Figure~\ref{f1}.

Taking the negative orientation for the surface $\partial U$
(i.e., choosing a unit normal vector $\nu$ on $\partial U$ that points inward), we can compute the flux of $fw$ across $\partial U$ as
$$
\int_{\partial U} fw\cdot \nu\, d\sigma=\int_{\mathcal D_0} fw\cdot \nu_0\,d\sigma - \int_{\mathcal D_1} fw\cdot \nu_1\,d\sigma\,,
$$
where $d\sigma$ denotes the induced surface measure and $\nu_s$ denotes
the unit normal on~$\cD_s$ pointing in the direction of $w$ (that
is, $w\cdot \nu_s>0$).

Using Eq.~\eqref{eqstep1}, the flux of $fw$ can also be written as
$$
\int_{\partial U} fw\cdot \nu\, d\sigma=\int_{\partial U} \curl
K(w)\cdot \nu\, d\sigma=0\,,
$$
with the integral vanishing by Stokes' theorem. Therefore we conclude
that the fluxes through the caps $\mathcal D_0$ and $\mathcal D_1$
must be equal, that is,
\begin{equation}\label{eqeqflux}
\int_{\mathcal D_{0}}fw\cdot \nu_0\, d\sigma=\int_{\mathcal D_{1}}fw\cdot \nu_1\, d\sigma\,.
\end{equation}

Suppose now that $f$ is not constant along the integral curves of
$w$. Then we can take a point $x_0\in D$ such that the function $f$
takes different values at the points $p_s:=\Phi^{-1}(s,x_0)\in \cD_s$, with
$s=0,1$. For concreteness, let us assume that 
\begin{equation}\label{eqfp}
f(p_0)<f(p_1)\,,
\end{equation}
the case $f(p_0)>f(p_1)$ being completely analogous. By the continuity
of $f$, we can then take the flow box narrow enough (i.e. with
$\mathcal D_0$ and $\mathcal D_1$ having very small diameters) such
that $c_{0}<c_{1}$, where 
$$c_{0}:=\max_{x\in \mathcal D_{0}} f(x)\,,\qquad c_{1}:=\min_{x\in
  \mathcal D_{1}} f(x)\,.
$$
Therefore, since $w\cdot \nu_s>0$ on $\cD_s$, we have the bound $$
\int_{\mathcal D_{0}}fw\cdot \nu_0\, d\sigma\leq c_{0}\int_{\mathcal D_{0}}w\cdot \nu_0\, d\sigma<c_{1}\int_{\mathcal D_{1}}w\cdot \nu_1\, d\sigma\leq\int_{\mathcal D_{1}}fw\cdot \nu_1\, d\sigma\,,
$$
where to obtain the second inequality we have used that, as $w$ is divergence-free, Stokes' theorem implies that
$$\int_{\mathcal D_{0}}w\cdot \nu_0\, d\sigma=\int_{\mathcal D_{1}}w\cdot \nu_1\, d\sigma\,.$$
This inequality above contradicts Eq.~\eqref{eqeqflux}, so we conclude
that $f$ must be constant along the integral curves of $w$, thus
proving that $f$ is a continuous first integral of~$w$ on $M\backslash
w^{-1}(0)$, as we had claimed. 

\subsubsection*{Step 3: There exists a continuous functional $\cC$ on $\Vect\backslash \{0\}$ such that derivatives of the invariant $\cI$ and of the helicity $\cH$ are related by
$
(D\cI)_w=\cC(w)(D\cH)_w
$
}

Let us start by noticing that Steps~1 and~2 imply that either $w$ has
a nontrivial first integral $F\in C^1(M)$ or the function $f$ defined
in Step~1 is a continuous first integral of $w$ in the complement of
its zero set. Now we observe that there exists a residual set $\cR$ of
vector fields in $\Vect$ such that any $w\in\cR$ is topologically
transitive and its zero set consists of finitely many hyperbolic
points. (We recall that a set is {\em residual}\/ if it is the intersection of
countably many open dense sets. In particular, a residual set is
always dense but not necessarily open.) This theorem was proved
in~\cite{Bessa} for divergence-free $C^1$ vector fields, not necessarily exact. However, it is not difficult to prove that the same result holds true for exact divergence free vector fields. Indeed, the proof of~\cite{Bessa} consists in perturbing a divergence-free vector field $w$ to obtain another divergence-free vector field $\tilde w$ of the form
$$
\tilde w=w+\sum_{i=1}^N v_i\,,
$$ 
where each $v_i$ is a $C^1$ divergence-free vector field supported in
a contractible set. Each vector field $v_i$ is necessarily exact
because any divergence-free vector field supported in a contractible
set is, so the resulting perturbed field $\tilde w$ is exact too. With
this observation, the main theorem in~\cite{Bessa} automatically
applies to the class of exact divergence-free $C^1$ vector fields, $\Vect$.

Hence let us take a vector field $w\in\cR$. Since it is topologically
transitive, it has an integral curve that is dense in $M$, so any
continuous first integral of $w$ must be a constant. Accordingly,
Steps~1 and~2 imply that $\curl K(w)=fw$ in $M\backslash w^{-1}(0)$,
with $f$ a first integral of $w$, and therefore the function $f$ is a
constant $c_w$ (depending on $w$) in the complement of the zero set
$w^{-1}(0)$. Since this set consists of finitely many points, $c_w$ is
the unique continuous extension of $f$ to the whole manifold $M$. As
$\curl K(w)$ is a continuous vector field,  for any $w\in\cR$ it follows that
\begin{equation}\label{curlKcw}
\curl K(w)=c_ww
\end{equation}
in $M$, so $\curl K(w)\times w=0$. 

Since the kernel $K$ is a continuous map $\Vect\to\Vect$, the fact
that $\curl K(w)\times w=0$ for all $w$ in the residual set
$\cR\subset\Vect$ implies that $\curl K(w)\times w=0$ for all
$w\in\Vect$. Therefore for any $w\in\Vect\backslash\{0\}$ we can define a function
$f\in C^0(M\backslash w^{-1}(0))$ by setting
\[
f:=\frac{w\cdot \curl K(w)}{|w|^2}\,,
\]
such that
\[
\curl K(w)=fw
\]
on $M\backslash w^{-1}(0)$. In view of the expression for $f$, the mapping $w\to f$ is continuous on $\Vect\backslash\{0\}$ due to the continuity of the kernel $K:\Vect\to\Vect$. Since $f$ is given by a $w$-dependent constant $c_w$ whenever $w$ lies in the residual set $\cR$ of $\Vect$,  we conclude that this must also be the case for all $w\in \Vect\backslash\{0\}$, so the
map $w\mapsto \frac12c_w$ defines a continuous functional
$\cC:\Vect\backslash\{0\}\to\RR$. (The factor $\frac12$ has been included for future
notational convenience.) The continuous functionals $\curl K(w)$ and $2\cC(w)\,w$ coinciding in a residual set, it stems that
for any $w\in\Vect\backslash\{0\}$ one has
\[
\curl K(w)=2\,\cC(w)\,w
\]
in all $M$.

Since the curl operator is invertible on $\Vect$ and $\cC(w)$ is just a
constant, we can use the above equation for $\curl K(w)$ to write the
derivative of $\cI$ at~$w$ as
$$
(D\cI)_w(u)=2\,\cC(w)\int_{M}\curl^{-1}w\cdot u\,.
$$
The claim of this step then follows upon recalling that the differential of the helicity is given by
$$
(D\cH)_w(u)=2\int_{M} \curl^{-1}w\cdot u\,.
$$

\subsubsection*{Step 4: The level sets of the helicity, $\cH^{-1}(c)$,
  are path connected subsets of $\Vect$.}

Let $w_0$ and $w_1$ be two vector fields in $\Vect$ with the same
helicity: 
\[
\cH(w_0)=\cH(w_1)=c\,.
\]
For concreteness, let us assume that $c$ is positive. It is easy to
see that the path connectedness of the level set $\cH^{-1}(c)$
is immediate if one can prove the existence of a path of positive helicity
connecting $w_0$ and~$w_1$, i.e., a continuous map
$w:[0,1]\to\Vect$ such that $w(0)=w_0$, $w(1)=w_1$ and $\cH(w(t))>0$
for all $t\in[0,1]$. Indeed, one can then set 
$$
\tilde w(t):=\bigg(\frac{c}{\cH(w(t))}\bigg)^{\frac12} \, w(t)
$$
to conclude that $\tilde w:[0,1]\to\Vect$ is a continuous path connecting $w_0$ and
$w_1$ of helicity~$c$: $\tilde w(0)=w_0$, $\tilde w (1)=w_1$ and $\cH(\tilde w (t))=c$ for all $t\in[0,1]$.

To show the existence of a path of positive helicity connecting
$w_0$ and $w_1$, we first observe that the curl defines a self-adjoint
operator with dense domain on the space of
exact divergence-free $L^2$ fields (see e.g.~\cite{GY}), so we can
take an orthonormal basis of eigenfields
$\{v_n^+,v_n^-\}_{n=1}^\infty$ satisfying $\curl v^{\pm}_n=\la_n^\pm
v_n^\pm$. Here we are denoting by $\la_n^+$ and $\la_n^-$ the positive
and negative eigenvalues of the curl, respectively. 

Given any vector field $v\in\Vect$,
we can expand $v$ in this orthonormal basis as
\[
v=\sum_{n=1}^\infty(c_n^+v_n^++c_n^-v_n^-)\,.
\]
This series converges in the Sobolev space $H^1$. As $\curl^{-1}
v_n^\pm=v_n^\pm/\la_n^\pm$, the helicity of the field $v$ can be written
in terms of the coefficients of the series expansion as
\begin{equation}\label{cHwj}
\cH(v)=\sum_{n=1}^\infty \Bigg(\frac{(c_n^+)^2}{\la_n^+}-\frac{(c_n^-)^2}{|\la_n^-|}\Bigg)\,. 
\end{equation}

We shall denote by $c_{j,n}^\pm$ the coefficients of the eigenfunction
expansion corresponding to $w_j$, with $j=0,1$. Let us fix two integers $n_j$ for which the
coefficient $c^+_{j,n_j}$ is nonzero (notice that the coefficients
corresponding to positive eigenvalues
cannot be all zero because of the formula~\eqref{cHwj} for the
helicity, which is positive in the case of $w_j$).

We can now construct the desired continuous path
$w:[0,1]\to\Vect$ of positive helicity connecting $w_0$ and $w_1$ by setting
\[
w(t):=\begin{cases}
8t\,c_{0,n_0}^+v_{n_0}^++(1-4t)\,  w_0 & \text{if } 0\leq t\leq
\frac14\,,\\[1mm]
2\cos(\pi t-\frac\pi4)\, c_{0,n_0}^+v_{n_0}^+ + 2\sin(\pi t-\frac\pi4)\, c_{1,n_1}^+v_{n_1}^+& \text{if } \frac14\leq t\leq
\frac34\,,\\[1mm]
(8-8t)\,c_{1,n_1}^+v_{n_1}^++(4t-3)\,  w_1 & \text{if } \frac34\leq t\leq 1\,.
\end{cases}
\]
Notice that $w(t)\in\Vect$ for all $t$ because both $w_j$ and the
eigenfields $v^+_{n_j}$
are in $\Vect$ (recall that the eigenfields of curl are automatically smooth
because they are also eigenfields of the Hodge Laplacian acting on
vector fields). It is also obvious that $w(0)= w_0$ and
$w(1)= w_1$. Furthermore, one can see that $w$ is a path
of positive helicity. For this, it is enough to use the
formula~\eqref{cHwj} for the helicity in terms of the coefficients of
the eigenfunction expansion. Indeed, since $\cH(w_j)=c$, an elementary computation then yields
\[
\cH(w(t))=\begin{cases}
16t\frac{(c_{0,n_0}^+)^2}{\la_0^+}+(1-4t)^2c & \text{if } 0\leq t\leq
\frac14\,,\\[1mm]
\frac{4(c_{0,n_0}^+)^2}{\la_{n_0}^+}\cos^2(\pi t-\frac\pi4)+ \frac{4(c_{1,n_1}^+)^2}{\la_{n_1}^+}\sin^2(\pi t-\frac\pi4)& \text{if } \frac14\leq t\leq
\frac34\,,\\[1mm]
16(1-t)\frac{(c_{1,n_1}^+)^2}{\la_{n_1}^+}+(4t-3)^2c & \text{if } \frac34\leq t\leq 1\,,
\end{cases}
\]
provided that $n_0\neq n_1$, so $\cH(w(t))>0$. When $n_0=n_1$, the only change in the formula above is that the value of $\cH(w(t))$ is
\[
\frac{4\Big(\cos(\pi t-\frac\pi4)c_{0,n_0}^++\sin(\pi t-\frac\pi4)c_{1,n_1}^+\Big)^2}{\la^+_{n_0}}
\]
if $\frac14\leq t\leq \frac34$, which is also positive. This proves the connectedness of $\cH^{-1}(c)$ when $c>0$.

The case where the constant~$c$ is negative is completely analogous
so, in order to finish the proof of the claim, it only remains to show
that the zero level set $\cH^{-1}(0)$ is path connected too. This is immediate
because two vector fields $w_0,w_1\in\Vect$ with $\cH(w_0)=\cH(w_1)=0$ can
be joined through the continuous path of zero helicity
$ w:[0,1]\to\Vect$ given by
\[
 w(t):=\begin{cases}
(1-2t)w_0 & \text{if } 0\leq t\leq
\frac12\,,\\
(2t-1)w_1 & \text{if } \frac12\leq t\leq1\,.
\end{cases}
\]
Obviously $ w(0)=w_0$, $ w(1)=w_1$ and $\cH( w(t))=0$ for all $t$, so the claim follows.

\subsubsection*{Step 5: The regular integral invariant $\cI$ is a function of the helicity.} 

We have shown in Step~3 that the derivatives of the functional $\cI$
and the helicity $\cH$ are related by $(D\cI)_w=\cC(w)(D\cH)_w$ at any
$w\in\Vect\backslash\{0\}$. In particular, this implies that $\cI$ is constant on
each path connected component of the level set $\cH^{-1}(c)\backslash\{0\}$. If $c\neq 0$, since $0$ is not contained in $\cH^{-1}(c)$, the aforementioned level set is path connected as proved in Step~4. The level set $\cH^{-1}(0)$ of zero helicity contains the $0$ vector field, so the set $\cH^{-1}(0)\backslash \{0\}$ does not need to be connected. However, since any component of $\cH^{-1}(0)\backslash \{0\}$ is path connected with $0$ as shown in the last paragraph of Step~4, the continuity of the functional $\cI$ in $\Vect$ implies that it takes the same constant value on any connected component of $\cH^{-1}(0)\backslash \{0\}$, so it is constant on the path connected level set $\cH^{-1}(0)$. We conclude that
there exists a function $f:\RR\to\RR$ which assigns a value of $\cI$
to each value of the helicity, i.e., $\cI=f(\cH)$. Moreover, $f$ is of
class $C^1$ because $\cI$ is a $C^1$ functional. The main theorem is
then proved.

\medskip

\begin{rremark}
  The only part of the proof where it is crucially used that the
  regularity of the vector fields is $C^1$ is in Step~3, when we
  invoke Bessa's theorem for generic vector fields in $\Vect$. To our
  best knowledge, it is not known if there is a residual subset of the
  space $\mathfrak X^k_{\rm ex}$ of exact divergence-free vector fields of class $C^k$,
  with $1<k\leq \infty$, whose elements do not
  admit a $C^{k-1}$ first integral. In particular, for $k>3$ the KAM theorem~\cite{KKP14} implies that there is no a residual subset of $\mathfrak X^k_{\rm ex}$ whose elements are topologically transitive vector fields, thus showing that Bessa's theorem does not hold for these spaces and hence it cannot be used to address the problem of the existence of a first integral for a generic vector field. Apart from the topological transitivity, we are not aware of other properties of a dynamical system implying that a vector field does not admit a (nontrivial) continuous first integral. The lack of results in this direction prevents us from extending the main theorem to regular
  integral invariants acting on $\mathfrak X^k_{\rm ex}$ with $k>1$.
\end{rremark}

\section*{Acknowledgments}

The authors are indebted to E. Kudryavtseva for sending us her paper~\cite{Ku} and explaining the main ideas of the proof of her theorem. We are also grateful to M. Bessa for his enlightening comments regarding the proof in~\cite{Bessa}. The authors are supported by the ERC Starting Grants~633152 (A.E.) and~335079
(D.P.-S.\ and F.T.L.) and by a fellowship from Residencia de Estudiantes (F.T.L.). This work is supported in part by the
ICMAT--Severo Ochoa grant
SEV-2011-0087.

\bibliographystyle{amsplain}

\end{document}